\input amstex
\documentstyle{amsppt}
\magnification=\magstep1                        
\hsize6.5truein\vsize8.9truein                  
\NoRunningHeads
\loadeusm

\magnification=\magstep1                        
\hsize6.5truein\vsize8.9truein                  
\NoRunningHeads
\loadeusm

\document
\topmatter

\title
On an Erd\H os Problem about the Maximum Modulus of Littlewood Polynomials on the Unit Circle
\endtitle

\rightheadtext{ultraflat unimodular polynomials}

\author Tam\'as Erd\'elyi
\endauthor

\address Department of Mathematics, Texas A\&M University,
College Station, Texas 77843, College Station, Texas 77843 \endaddress

\endthanks


\date May 12, 2026
\enddate

\email terdelyi\@tamu.edu
\endemail


\abstract
Making a progress in an Erd\H os problem dating back to 1957, we show that 
$$\max_{t \in {\Bbb R}}{|P_n(e^{it})|^2} \geq n+1 + \frac{1}{38}n^{1/3}$$
for every polynomial $P_n$ of degree $n$ with all of its coefficients in $\{-1,1\}$.
\endabstract
\endtopmatter

\head 1. Introduction \endhead
Let
$${\Cal K}_n := \left\{P_n: P_n(z) = \sum_{k=0}^n{a_k z^k}, \enskip  a_k \in {\Bbb C}\,,\enskip |a_k| = 1 \right\}\,.$$
The class ${\Cal K}_n$ is often called the collection of all (complex) unimodular polynomials of degree $n$.
Let
$${\Cal L}_n := \left\{P_n: P_n(z) = \sum_{k=0}^n{a_k z^k}, \enskip  a_k \in \{-1,1\} \right\}\,.$$
The class ${\Cal L}_n$ is often called the collection of all (real) unimodular polynomials of degree $n$.
The elements of the class ${\Cal L}_n$ are also called Littlewood polynomials of degree $n$.
By Parseval's formula,
$$\int_{0}^{2\pi}{\left| P_n(e^{it}) \right|^2 \, dt} = 2\pi(n+1)$$
for all $P_n \in {\Cal K}_n$. Therefore
$$\min_{t \in {\Bbb R}}{|P_n(e^{it})|} \leq \sqrt{n+1} \leq \max_{t \in \Bbb R}{|P_n(e^{it})|} \tag 1.1$$
for all $P_n \in {\Cal K}_n$. An old problem (or rather an old theme) is the following.

\proclaim{Problem 1.1 (Littlewood's Flatness Problem)}
How close can a polynomial $P_n \in {\Cal K}_n$ or $P_n \in {\Cal L}_n$ come to satisfying
$$|P_n(e^{it})| = \sqrt{n+1}\,, \qquad t \in {\Bbb R}? \tag 1.2$$
\endproclaim

Obviously (1.1) is impossible if $n \geq 1$. So one must look for less than (1.2), but then there 
are various ways of seeking such an ``approximate situation". One way is the following.
In his paper [Li66] Littlewood had suggested that, conceivably, there might exist a sequence
$(P_n)$ of polynomials $P_n \in {\Cal K}_n$ (possibly even $P_n \in {\Cal L}_n$) such
that $(n+1)^{-1/2}|P_n(e^{it})|$ converge to $1$ uniformly in $t \in {\Bbb R}$.
We shall call such sequences of unimodular polynomials ``ultraflat". More precisely, we give the
following definition.

\proclaim{Definition 1.2} Given a positive number $\varepsilon$, we say that a polynomial 
$P_n \in {\Cal K}_n$ is $\varepsilon$-flat if
$$(1-\varepsilon)\sqrt{n+1} \leq |P_n(e^{it})| \leq (1 + \varepsilon)\sqrt{n + 1}\,,
\qquad t \in {\Bbb R}\,.$$
\endproclaim

\proclaim{Definition 1.3} Let $(n_k)$ be an increasing sequence of positive integers. Given a sequence 
$(\varepsilon_{n_k})$ of positive numbers tending to $0$, we say that a sequence $(P_{n_k})$ of polynomials 
$P_{n_k} \in {\Cal K}_{n_k}$ is $(\varepsilon_{n_k})$-ultraflat if each $P_{n_k}$ is $\varepsilon_{n_k}$-flat. 
We simply say that a sequence $(P_{n_k})$ of polynomials $P_{n_k} \in {\Cal K}_{n_k}$ is
ultraflat if it is $(\varepsilon_{n_k})$-ultraflat with a suitable sequence $(\varepsilon_{n_k})$ 
of positive numbers tending to $0$.
\endproclaim

The existence of an ultraflat sequence of unimodular polynomials seemed very unlikely, in view of a
1957 conjecture of P. Erd\H os (Problem 22 in [Er57]) asserting that, for all
$P_n \in {\Cal K}_n$ with $n \geq 1$,
$$\max_{t \in {\Bbb R}}{|P_n(e^{it})|} \geq (1 + \varepsilon) \sqrt{n+1}\,, \tag 1.3$$
where $\varepsilon > 0$ is an absolute constant (independent of $n$).
Yet, refining a method of K\"orner [K\"o80], Kahane [Ka80] proved that there exists
a sequence $(P_n)$ with $P_n \in {\Cal K}_n$ which is $(\varepsilon_n)$-ultraflat, where
$\varepsilon_n = O\left(n^{-1/17} \sqrt{\log n} \right)\,.$
(Kahane's paper contained though a slight error which was corrected in [QS96].)  
Thus the Erd\H os conjecture (1.3) was disproved for the classes ${\Cal K}_n$.
For the more restricted class ${\Cal L}_n$ the analogous Erd\H os conjecture
is unsettled to this date. It is a common belief that the analogous Erd\H os conjecture
for ${\Cal L}_n$ is true, and consequently there is no ultraflat sequence of polynomials
$P_n \in {\Cal L}_n$.
An interesting result related to Kahane's breakthrough is given in [Be91].
For an account of some of the work done till the mid 1960's, see Littlewood's book [Li68]
and [QS96]. Littlewood polynomials with small $L_4$ norm have also 
been intensively studied, see [BM00], [BC00], [BC01],[BC02], and [JK13], for example.   

The structure of ultraflat sequences of unimodular polynomials is studied in
[Er00a], [Er00b], [Er01a], [Er01b], [Er02], [Er03], and [Er21] where several 
conjectures of Saffari [Sa92] and Queffelec and Saffari [QS96] are proved. 
A recent paper of Bombieri and Bourgain [BB09] is devoted to the construction of ultraflat 
sequences of unimodular polynomials. In particular, they obtained a much improved estimate 
for the error term. A major part of their paper deals also with the long-standing problem 
of the effective construction of ultraflat sequences of unimodular polynomials.  
A recent breakthrough result by Balister, Bollob\'as, Morris, 
Sahasrabudhe, and Tiba [BB20] states that flat Littlewood polynomials exist, more precisely there exist 
absolute constants $\eta_2 > \eta_1 > 0$ and a sequence $(P_n)$ of Littlewood polynomials $P_n$ 
of degree $n$ such that
$$\eta_1 \sqrt{n} \leq |P_n(e^{it})| \leq \eta_2 \sqrt{n}\,, \qquad t \in {\Bbb R}\,,$$
confirming a conjecture of Littlewood [Li66] dating back to 1966. Moreover, it is shown in [Er22] that the  
sequence $(P_n)$ of Littlewood polynomials $P_n$ of degree $n$ can be chosen in a way that in addition to 
the above flatness properties a certain symmetry is satisfied by the coefficients of $P_n$ making the Littlewood 
polynomials $P_n$ close to skew-reciprocal.

There are quite a few recent publications on or related to ultraflat sequences of unimodular 
polynomials. Some of them (not mentioned before) are are [Bo02], [Sa01], [QS95], [Od18], and [Mo]. 

\head 2. New Results \endhead

\proclaim{Theorem 2.1}
We have  
$$\max_{t \in {\Bbb R}}{\left|P_n(e^{it})\right|^2} \geq n+1 + \frac{1}{38}n^{1/3}$$
for every polynomial $P_n \in {\Cal L}_n$.
\endproclaim

\proclaim{Theorem 2.2}
We have 
$$\max \left( \max_{t \in {\Bbb R}}{\left|P_n(e^{it})\right|^2} - (n+1), (n+1) - \min_{t \in {\Bbb R}}{\left|P_n(e^{it})\right|^2} \right) \geq \left(n/2\right)^{1/2}$$
for every $P_n \in {\Cal L}_n$.
\endproclaim

\proclaim{Theorem 2.3}
We have
$$\min_{t \in {\Bbb R}}{\left|P_n(e^{it}\right)|^2} \leq (n+1) - \frac{\log n}{4\pi} < (n+1) + \frac{\log n}{4\pi}  \leq \max_{t \in {\Bbb R}}{\left|P_n(e^{it})\right|^2}$$
for every polynomial $P_n \in {\Cal L}_n$.
\endproclaim 

Theorem 2.3 may be compared with the Appendix of [Bo93] and Theorem 1.5 in [Er11].

A polynomial $P_n$ of the form
$$P_n(z) = \sum_{k=0}^n {a_k z^k}\,, \qquad a_k \in {\Bbb C}\,, \tag 2.1$$
is called conjugate-reciprocal if 
$$a_{n-k} = \overline{a}_k\,, \qquad k=0,1,\ldots,n\,. \tag 2.2$$ 

The proof of the remark below may be found in [BE03], but in Section 4 we present
a short self-contained proof of it. 
 
\proclaim{Remark 2.4}
We have 
$$\max_{t \in {\Bbb R}}{\left|P_n(e^{it})\right|} \geq (4/3)^{1/2} n^{1/2}$$ 
for every conjugate-reciprocal polynomial $P_n \in {\Cal K}_n$.
\endproclaim

In fact, in [Er15], Theorem 2.6] we proved the result below.

\proclaim{Remark 2.5} 
There is an absolute constant $\varepsilon > 0$ succh that
$$\max_{t \in {\Bbb R}}{\left|P_n(e^{it})\right|} \geq (1 + \varepsilon) \left(4/3\right)^{1/2} n^{1/2}$$
for every conjugate-reciprocal polynomial $P_n \in {\Cal K}_n$ and for all sufficiently large $n$.
\endproclaim

The proof of the remark below follows simply from well known results, but in Section 4 we present 
a short self-contained proof of it.  

\proclaim{Remark 2.6}
We have
$$0 = \min_{t \in {\Bbb R}}|P_n(e^{it})|$$
for every conjugate-reciprocal $P_n \in {\Cal K}_n$, $n \geq 1$.
\endproclaim

\head 3. Lemmas \endhead

To prove Theorem 2.1 we need the Bernstein inequality in $L_q$ for trigonometric polynomials in $L_q$, $q > 0$.
In fact, we need only the case $q=1$ of it. See [Ar79] and [Ar81] for the general case $q > 0$. For a book proof see 
[DL93] or [Er20]. For the case $q \geq 1$ see [BE95]. 

\proclaim{Lemma 3.1} We have 
$$\int_{0}^{2\pi}{\left|Q_n^{\prime}(t)\right|^q \, dt} \leq n^q \int_{0}^{2\pi}{\left|Q_n(t)\right|^q \, dt}$$
for every trigonometric polynomial $Q_n$ of degree at most $n$ with complex coefficients.
\endproclaim

We call the functions 
$$T_n(t) = a_0 + \sum_{k=1}^n{a_k \cos(kt) + b_k \sin(kt)}\,, \qquad a_k, b_k \in {\Bbb R}, \quad a_nb_n \neq 0\,,$$
a real trigonometric polynomial of degree $n$. The following inequality is a version of the Bernstein-Szeg\H o inequality.
For a proof of it see [BE95]. We will use it in the proof of Theorem 2.1 as well.

\proclaim{Lemma 3.2}
Suppose $T_n$ is a real trigonometric polynomial of degree at most $n$ satisfying $m \leq T_n(t) \leq M$ for all $t \in {\Bbb R}$.
We have 
$$\left|T_n^{\prime}(t)\right|^2 + n^2\left|T_n(t) - \frac{M+m}{2}\right|^2 \leq n^2\left(\frac{M-m}{2}\right)^2\,, \qquad t \in {\Bbb R}\,.$$
\endproclaim

A short proof of Remark 2.3 is based on the following Bernstein-type inequality, proved in [BE95], for conjugate-reciprocal polynomials.

\proclaim{Lemma 3.3} 
We have 
$$\max_{t \in {\Bbb R}}{\left|P_n^{\prime}(e^{it})\right|} \leq \frac n2 \max_{t \in {\Bbb R}}{\left|P_n(e^{it})\right|}$$
for every conjugate-reciprocal polynomials $P_n$ of degree $n$.
\endproclaim

To prove Theorem 2.3 we need Pavlovi\'c's improvement of Hardy's inequality below. For a proof See Theorem 2.15 in [Pa14].

\proclaim{Lemma 3.4}
We have
$$\sum_{k=0}^m{\frac{|a_k|}{k+1}} \leq \int_{0}^{2\pi}{|Q_m(e^{it})| \,dt}$$
for every polynomial $Q_m$ of the form 
$$Q_m(z) = \sum_{k=0}^m{a_k z^k}\,, \qquad a_k \in {\Bbb C}\,.$$
\endproclaim

\proclaim{Corollary 3.5}
We have
$$\frac 12 \log n \leq \int_{0}^{2\pi}{|Q_{2n}(e^{it})| \,dt}$$
for every polynomial $Q_{2n}$ of the form
$$Q_{2n}(z) = \sum_{k=0}^{2n}{b_k z^k}\,, \qquad a_k \in {\Bbb C}\,, \quad |a_{2u}| \geq 1\,, \quad u=0,1,\ldots,n\,, \enskip 2u \neq n\,. \tag 3.1$$

\endproclaim

\head 4. Proofs \endhead

\demo{Proof of Theorem 2.1}
Let $P_n \in {\Cal L}_n$. Let $T_n$ be the trigonometric polynomial of degree $n$ defined by
$$T_n(t) := \left|P_n(e^{it})\right|^2 - (n+1) = P_n(e^{it})P_n(e^{-it}) - (n+1)\,, \qquad t \in {\Bbb R}\,. \tag 4.1$$
We have 
$$T_n(t) =\sum_{k=0}^{2n}{b_{k-n}e^{i(k-n)kt}}\,, \quad t \in {\Bbb R}\,.$$
Observe that if $0 \leq k \neq 2n$ and  $k$ is an even integer, then $b_{k-n}$ is a sum of numbers in $\{-1,1\}$ with an odd number of terms, 
hence $b_{k-n}$ is an odd integer, so $|b_{k-n}| \geq 1$. Therefore by Parseval's formula we get that
$$\int_{0}^{2\pi}{\left|T_n^{\prime}(t)\right|^2 \, dt} = 2\pi \sum_{k=0}^{2n}{|(k-n)b_{k-n}|^2} \geq 2\pi \frac{n(n+1)(2n+1)}{6}\,. \tag 4.2$$  
We have
$$\int_E{\left|T_n^{\prime}(t)\right|^2 \, dt} \leq \left( \int_{0}^{2\pi}{|T_n^{\prime}(t)| \, dt} \right) \left( \max_{t \in E}{|T_n^{\prime}(t)|} \right) \tag 4.3$$ 
for every measurable set $E \subset {\Bbb R}$.
Suppose that there are $0 \leq \delta_n \leq (n+1)/16$ such that 
$$\max_{t \in {\Bbb R}}{\left|P_n(e^{it})\right|^2} \leq n+1 + \delta_n\,. \tag 4.4$$ 
As
$$\int_{0}^{2\pi}{T_n(t) \, dt} = 0$$
we have
$$\int_{0}^{2\pi}{T_n^-(t) \, dt} = \int_{0}^{2\pi}{T_n^+(t) \, dt}\,,$$
where $T_n^-(t) := \max(-T_n(t),0)$ and $T_n^+(t) := \max(T_n(t),0)\,.$
Observe that (4.4) implies that
$$0 \leq \int_{0}^{2\pi}{T_n^-(t) \, dt} = \int_{0}^{2\pi}{T_n^+(t) \, dt} \leq 2\pi \delta_n\,. \tag 4.5$$

Using the Bernstein inequality in $L_1$ for trigonometric polynomials $Q_n$ of degree $n$ defined by 
$$Q_n(t) :=  n+1 + \delta_n - P_n(e^{it})P_n(e^{-it}) = n+1 + \delta_n - \left|P_n(e^{it})\right|^2 \geq 0\,, \qquad t \in {\Bbb R}\,, \tag 4.6$$
we have
$$\int_{0}^{2\pi}{|Q_n^{\prime}(t)| \, dt} \leq n \int_0^{2\pi}{|Q_n(t)| \, dt}\,.$$
Hence, recalling (4.1), (4.5), and (4.6),  we get
$$\int_{0}^{2\pi}{|T_n^{\prime}(t)| \, dt} \leq 2\pi n\delta_n\,. \tag 4.7$$
Let $k \geq 1$ be the largest integer such that $2^k \delta_n < (n+1+\delta_n)/2$. We define the sets
$$\eqalign{B_1 := & \{t \in [0,2\pi]: - \delta_n <  T_n(t) \leq \delta_n\}\,, \cr 
           A_j := & \{t \in [0,2\pi]: - 2^j\delta_n <  T_n(t) \leq - 2^{j-1}\delta_n\}\,, \qquad j=1,2,\ldots,k\,, \cr 
           B 2 := & \{t \in [0,2\pi]: -(n+1) \leq T_n(t) \leq -(n+1)/4\}\,. \cr}$$
Observe that 
$$[0,2\pi] =  B_1 \cup B_2 \cup_{j=1}^k{A_j}\,, \tag 4.8$$
$$m(B_1) \leq 2\pi\,, \tag 4.9$$
$$m(A_j)2^{j-1}\delta_n \leq \int_{A_j}{-T_n(t) \, dt} \leq \int_0^{2\pi}{T_n^-(t) \, dt} \leq 2\pi \delta_n \,,$$ 
and hence
$$m(A_j) \leq 4\pi 2^{-j}\,, \qquad j=1,2,\ldots,k\,, \tag 4.10$$
where $m(E)$ denotes the linear Lebesgue measure of the set $E \subset {\Bbb R}$.
Observe also that
$$m(B_2)(n+1)/4 \leq \int_{B_2}{-T_n(t) \, dt} \leq \int_{0}^{2\pi}{T_n^-(t) \, dt} \leq 2\pi \delta_n \,,$$
and hence
$$m(B_2) \leq 8\pi \frac{\delta_n}{n+1}\,. \tag 4.11$$

Applying the Bernstein-Szeg\H o inequality (Lemma 3.2) with $T_n$, $m := -(n+1)$, and $M := \delta_n$, and recalling (4.4), we get
$$\left|T_n^{\prime}(t)\right|^2 + n^2\left|T_n(t) - \frac{M+m}{2}\right|^2 \leq n^2\left(\frac{M-m}{2}\right)^2\,, \qquad t \in {\Bbb R}\,, \tag 4.12$$

If $t \in B_1$, then (4.12) implies
$$\left|T_n^{\prime}(t)\right|^2 + n^2\left( -\delta_n - (\delta_n - (n+1))/2 \right)^2 \leq n^2\left((\delta_n + (n+1))/2 \right)^2\,.$$
Hence for $t \in B_1$ we have
$$\split \left|T_n^{\prime}(t)\right|^2 \leq & n^2\left( \left((\delta_n/2 + (n+1))/2 \right)^2 - n^2\left( (n+1)/2 - 3\delta_n/2 \right)^2 \right) \cr
\leq & n^2((n+1) - \delta_n)2\delta_n = 2n^2(n+1)\delta_n\,, \cr \endsplit$$
from which
$$\left|T_n^{\prime}(t)\right| \leq 2^{1/2}n(n+1)^{1/2}\delta_n^{1/2}\,, \qquad t \in B_1\,, \tag 4.13$$
follows.

If $t \in A_j, j=1,2,\ldots,k$, then (4.12) implies
$$\left|T_n^{\prime}(t)\right|^2 + n^2\left(-2^j \delta_n - (\delta_n - (n+1))/2 \right)^2 \leq n^2\left((\delta_n + (n+1))/2 \right)^2\,.$$
Hence for $t \in A_j, j=1,2,\ldots,k$, we have
$$\split \left|T_n^{\prime}(t)\right|^2 \leq & n^2\left( \left(\delta_n/2 + (n+1)/2 \right)^2 - \left(-2^j\delta_n - (\delta_n - (n+1))/2 \right)^2 \right) \cr 
= & n^2 \left( (n+1) - 2^j\delta_n \right)\left( \delta_n + 2^j\delta_n \right)\,, \cr \endsplit$$
from which
$$\left|T_n^{\prime}(t)\right| \leq 2^{1/2}n(n+1)^{1/2}2^{j/2}\delta_n^{1/2}\,, \qquad t \in A_j\,, \quad j=1,2,\ldots,k\,,\tag 4.14$$
follows.  

If $t \in B_2$, then (4.12) implies
$$\left|T_n^{\prime}(t)\right|^2 \leq n^2\left(((n+1) + \delta_n)/2 \right)^2\,,$$
from which
$$\left|T_n^{\prime}(t)\right| \leq n(n+1)\,, \qquad t \in B_2\,. \tag 4.15$$

Using (4.2), (4.3), (4.7)--(4.15), we get

$$\split & \frac{2\pi n(n+1)(2n+1)}{6} \leq \int_{0}^{2\pi}{\left|T_n^{\prime}(t)\right|^2 \, dt} \cr
\leq & \int_{B_1}{\left|T_n^{\prime}(t)\right|^2 \, dt} + \sum_{j=1}^k{\int_{A_j}{\left|T_n^{\prime}(t)\right|^2 \, dt}} + \int_{B_2}{\left|T_n^{\prime}(t)\right|^2 \, dt} \cr 
\leq & \left( \int_{0}^{2\pi}{|T_n^{\prime}(t)| \, dt} \right) \cdot \cr  
& \left(m(B_1) \max_{t \in B_1}{|T_n(t)|} + \sum_{j=1}^k{m(A_j) \max_{t \in A_j}{\left|T_n(t)\right|}} + m(B_2) \max_{t \in B_2}{|T_n(t)|} \right) \cr  
\leq & 2\pi n\delta_n \left( 2\pi 2^{1/2}n(n+1)^{1/2}\delta_n^{1/2} + \sum_{j=1}^k{4\pi 2^{-j} 2^{1/2}n(n+1)^{1/2}2^{j/2}\delta_n^{1/2}} + 8\pi \frac{\delta_nn(n+1)}{n+1} \right) \cr 
= & 2\pi n^2(n+1)^{1/2}\delta_n^{3/2} \left( 2\pi 2^{1/2} + 4\pi 2^{1/2}\sum_{j=1}^k{2^{-j/2}} + 8\pi \delta_n^{1/2}(n+1)^{-1/2} \right) \cr 
\leq & 2\pi n^2(n+1)^{1/2}\delta_n^{3/2} 2\pi \left( 1.42 + 6.84 + 4 \right)\,, \cr \endsplit$$
hence, choosing $\delta_n := cn^{1/3}$, we get 
$$\frac 13 < c^{3/2} \pi 24.52\,.$$
Thus $\displaystyle{\left( \frac{1}{\pi 73.56} \right)^{2/3}} < c$, which is impossible if $\displaystyle{c := \frac{1}{38} < \left( \frac{1}{\pi 73.56} \right)^{2/3}}$.
\qed \enddemo

\demo{Proof of Theorem 2.2}
Let $P_n \in {\Cal L}_n$. As in the proof of Theorem 2.1 we have
$$\left|P_n(e^{it})\right|^2 - (n+1) = P_n(e^{it})P_n(e^{-it}) - (n+1) = \sum_{k=0}^{2n}{b_{k-n}e^{ikt}}\,, \quad t \in {\Bbb R}\,. $$
where if $0 \leq k \neq 2n$ and $k$ is an even integer, then $b_{k-n}$ is a sum of numbers in $\{-1,1\}$ with an odd number of terms, hence $b_k$ is an odd integer, 
so $|b_{k-n}| \geq 1$. Therefore by Parseval's formula we get that  
$$\int_{0}^{2\pi}{\left(\left|P_n(e^{it})\right|^2 - (n+1) \right)^2 \, dt} = 2\pi \sum_{k=0}^{2n}{|b_k|^2} \geq 2\pi n\,. \tag 4.16$$
Now let 
$$f_n(t) := \left|P_n(e^{it})\right|^2 - (n+1) \,,$$
$f_n^+(t) := \max(f_n(t),0)$ and $f_n^-(t) := \max(-f_n(t),0)\,.$
Observe that 
$$f_n(t)^2 := f_n^+(t)^2 + f_n^{-}(t)^2$$ 
and by (4.16) we have either $\displaystyle{\int_{0}^{2\pi}{f_n^+(t)^2 \, dt} \geq \pi n}$ or 
$\displaystyle{\int_{0}^{2\pi}{f_n^-(t)^2 \, dt} \geq \pi n}$, and the theorem follows. 
\qed \enddemo

\demo{Proof of Theorem 2.3}
As in the proof of Theorem 2.1 we have
$$f_n(t) := \left|P_n(e^{it})\right|^2 - (n+1) = P_n(e^{it})P_n(e^{-it}) - (n+1) = e^{-int}Q_{2n}(e^{it})\,, \quad t \in {\Bbb R}\,, \tag 4.17$$
where $Q_{2n}$ is a polynomial of the form (3.1). Observe that $\displaystyle{\int_{0}^{2\pi}{f_n(t) \, dt} = 0}$, hence using the notation 
introduced in the proof of Theorem 2.2, we have
$$\int_{0}^{2\pi}{f_n^+(t) \, dt} = \int_{0}^{2\pi}{f_n^-(t) \, dt} = \frac 12 \int_{0}^{2\pi}{|f_n(t)| \, dt}\,,$$
hence the theorem follows from (4.17) and Corollary 3.5. 
\qed \enddemo

\demo{Proof of Remark 2.4}
Let $P_n \in {\Cal K}_n$ be self-reciprocal. Combining Parseval's formula with Lemma 3.3, we have
$$\split 2\pi \frac{n^3}{3} \leq & 2\pi \frac{n(n+1)(2n+1)}{6} = 2\pi \sum_{k=1}^n{k^2} = \int_0^{2\pi}{\left|P_n^{\prime}(e^{it})\right|^2 \, dt} \cr 
& \leq 2\pi \left( \frac n2 \right)^2 \left( \max_{t \in {\Bbb R}}{\left|P_n(e^{it})\right|} \right)^2 \,, \cr \endsplit$$
and the result follows. 
\qed \enddemo

\demo{Proof of Remark 2.6}
Leq $n \geq 1$ be an integer and let $P_n \in {\Cal K}_n$ be conjugate-reciprocal of the form (2.1).  

First let $m=2m+1$ be odd. Let 
$$a_{m-j} := e^{-i\gamma_j}\,, \qquad \gamma_j \in [0,2\pi)\,, \quad j=0,1,\ldots,m\,.$$
Then 
$$a_{m+1+j} := e^{i\gamma_j}\,, \qquad \gamma_j \in [0,2\pi)\,, \quad j=0,1,\ldots,m\,.$$
We define
$$F(t) := e^{-(m+1/2)t}P_n(e^{it}) = \sum_{j=0}^m{2\cos((2j+1)t/2+\gamma_j})\,.$$
Observe that 
$$\int_0^{2\pi}{F(2t)} = 0\,,$$
which implies that $F(t)$ and hence $P_n(e^{it})$ vanishes at some $t_0 \in [0,4\pi)$.

Now let $n=2m$ be even. 
We define
$$e^{-mt}P_n(e^{it}) =: 2\Re(G(e^{it}))\,, \qquad t \in {\Bbb R}\,,$$ 
where
$$G(z) := 2a_m + \sum_{j=1}^{m}\left(a_{m+j}jz^j\right), \qquad z \in {\Bbb C}\,.$$
Observe that $G$ is a polynomial with constant term $2a_m$ and leading coefficient $2a_{2m}$, where $a_m \in \{-1,1\}2$ and $|a_{2m}| = 1$. 
Hence there is at least one zero of $G$ in the open unit disk of the complex plane. Therefore by the Argument Principle 
$G(e^{it})$ goes around $G(0) = 2a_m \in \{-2,2\}$ at least once as $t$ increases from $0$ to $\pi$. We conclude that  
$\Re(G(e^{it}))$, and hence $P_n(e^{it})$ as well, vanish in $[0,2\pi)$ at least twice. 
\qed \enddemo

\enddocument